\documentclass{article}
\usepackage{amsmath, amsfonts, amssymb}  
\title{Extreme and Exposed Points Arising from Rational Kernels}
\author{Steve Fisher, Northwestern University (retired)}
\begin{document}
\maketitle

$\Delta$ denotes the open unit disc $\{z: |z| < 1 \}$ in the complex plane and \textbf{T} its boundary, the unit circle \textbf{T} $= \{z : |z|=1 \}$. If $P$ is a polynomial of exact degree $n$, we define a related polynomial $\widetilde{P}$ by \begin{equation} \widetilde{P}(z) = z^n \overline{P(1/\overline{z})}. \end{equation}

\bigskip

\textbf{Definition} Let $\Lambda$ be a convex set in a linear space. A point $\vec P \in \Lambda$ is an \textit{extreme point} if it is not in the interior of a line segment in $\Lambda.$ That is, \begin{displaymath} \textrm{if   } \vec P = \lambda \vec P_1 + (1-\lambda) \vec P_2, 0 < \lambda < 1; \vec P_1, \vec P_2 \in \Lambda, \textrm{   then  } \vec P_1=\vec P_2. \end{displaymath} $\vec P$ is an \textit{exposed point} of $\Lambda$ if there is a tangent hyperplane to $\Lambda$ that intersects $\Lambda$ exactly at $\vec P$. Clearly an exposed point of $\Lambda$ is also an extreme point. The converse does not hold in general. 

\bigskip

In what follows we investigate the extreme and exposed points of the convex set in $\mathbb{C}^n$ arising from a rational kernel and the Hardy space $H^1$. Specifically, let $\beta_1,...,\beta_n$ be distinct points in $\Delta$, none of which is the origin, and define a closed convex set in $\mathbb{C}^{n}$ by  \begin{equation} \Lambda = \{ (f(\beta_1),...,f(\beta_n)): f \in H^1, ||f||_1 \le 1 \}. \end{equation} Duren [1; page 138] has a full exposition of those $H^1$ functions $F$ that produce boundary points of $\Lambda$ since $(F(\beta_1),...,F(\beta_n))$ is in the boundary of $\Lambda$ if and only if $F$ is a solution of an extremal problem of the form \begin{equation} \sup \{ \textrm{Re  } [\sum_{j=1}^n c_j f(\beta_j)] : f \in H^1, ||f||_1 \le 1 \} \end{equation} for an appropriate choice of complex scalars $c_1,..,c_n$. Evidently, $\Lambda$ contains a neighborhood of the origin in $\mathbb{C}^n$.

\bigskip

\textbf{N.B.} Proposition 1 below can be found, for instance, in [1]. We give a proof to establish notation and relationships. S. Ya. Havinson [2] has discussed the issue of extreme points on more general planar domains; some of his results are contained in Theorem 2. I am indebted to Dimitry Khavinson for this reference. Duren [1] does not use the term "exposed point" but his observations on page 139 about the uniqueness of the solution in the case $p=1$ seem to encompass what is in Theorem 4 and Corollary 5. Consequently, any virtue of what follows is due to the consistency of notation and not to the originality of the results. All integrals are taken over the unit circle \textbf{T} with respect to Lebesgue measure $\sigma.$

\bigskip

\textbf{Proposition 1} Suppose $F \in H^1, ||F||_1=1$. Then $(F(\beta_1),...,F(\beta_n))$ lies in the boundary of $\Lambda$ if and only if $F$ is of  the form \begin{equation} F(z) = z^{\ell}W(z)q^2(z) \prod_{j=1}^n (1 - \overline{\beta_j}z)^{-2} \end{equation} where $q$ is a polynomial of degree $n-1$ or less with no zeros in $\Delta$, $0 \le \ell \le n-\textrm{degree  }q$, and $W$ is a Blaschke product whose zeros, if any, are also zeros of $\widetilde{q}$ in $\Delta$.

\medskip

\textbf{Proof}  For simplicity of notation, we define \begin{equation} B(z) = \prod_{k=1}^n \frac{z- \beta_k}{1 - \overline{\beta_k}z} \textrm{   and   } P(z) =\prod_{k=1}^n (1 - \overline{\beta_k}z)^{-2}. \end{equation} We note for future use that \begin{equation} \frac{e^{int}P(e^{it})}{B(e^{it})} > 0 \textrm{   on \textbf{T}}. \end{equation} Suppose $(F(\beta_1),...,F(\beta_n))$ is a boundary point of $\Lambda$. Then there are complex scalars $c_1,...,c_n$, not all of which are zero, such that \begin{equation}  1= \sum_{k=1}^n c_k F(\beta_k) \ge \textrm{Re   }\sum_{k=1}^n f(\beta_k), \textrm{   for all  } f \in H^1, ||f||_1\le 1. \end{equation} Let $R$ be the rational function  \begin{equation} R(z) = z \sum_{k=1}^n c_k \frac{1}{z-\beta_k} \end{equation} and let $g_0$ be a best approximation to $R$ in $L^{\infty}$ from $H^{\infty}_0$; finally, set \begin{displaymath} K=R-g_0. \end{displaymath}  The function $K$ is called the \textit{kernel function} for $F$. Thus, $K$ has a pole of order 1 at the point $\beta_k, k=1,..,n$ (unless $c_k=0$.) Moreover, $K(z)$ is bounded and analytic in the ring $1-\delta \le |z| < 1$ for any 
sufficiently small $\delta$. The inequality expressed in (7) may be rewritten as \begin{displaymath} \int FK d\sigma \ge \textrm{Re  } \int Kf d \sigma,  \textrm{   for all  } f \in H^1, ||f||_1\le 1. \end{displaymath}Duality then gives  \begin{displaymath} FK \ge 0  \textrm{   and  } |K|=1 \textrm{   a.e. on  } \textbf{T}. \end{displaymath}Therefore, \begin{displaymath} F(e^{it})K(e^{it}) e^{-int} \prod_{k=1}^n (e^{it}-\beta_k)(1- \overline{\beta_k} e^{it}) \ge 0 \textrm{   a.e. on  } \textbf{T}. \end{displaymath} However, the function \begin{displaymath}F(z)K(z)  \prod_{k=1}^n (z-\beta_k)(1- \overline{\beta_k} z) \end{displaymath} lies in $H_0^1$. Let $m \ge 1$ be the order of the zero of this function at the origin. It follows from the Fejer-Riesz Theorem that there is a polynomial $q$ of 
degree $n-m$ with no zero in $\Delta$ such that  \begin{displaymath} F(z)K(z) \prod_{j=1}^n (z- \beta_j)(1-\overline{\beta_j} z) = z^m q(z) \widetilde{q}(z). \end{displaymath}  From (5), (6), and the above we obtain  \begin{equation} F(z)K(z)B(z) = z^m q(z) \widetilde{q}(z) P(z) = z^m \frac{\widetilde{q}(z)}{q(z)} q^2(z)P(z).  \end{equation} Let $F=IG$ be the inner-outer factorization of $F$. Since $KB$ is inner and $q^2P$ is outer, when we equate the inner and outer factors from the far sides of (9) we obtain \begin{equation} I(z)K(z)B(z)= z^m \frac{\widetilde{q}(z)}{q(z)} \textrm{   and   } G(z) = q^2(z)P(z). \end{equation} Hence, $I$ and $KB$ are both finite Blaschke products. We define $\ell$ to be the order of the zero of $I$ at the origin; any remaining zeros of $I$ must be zeros of $\frac{\widetilde{q}(z)}{q(z)}$. Thus, $F$ has the form given in (4).

\medskip

Conversely, suppose $F$ has unit norm and is of the form in (4):  \begin{displaymath} F(z) = z^{\ell} W(z) q^2(z)P(z). \end{displaymath} If $q$ has degree $n-m$, then $\ell \le n - \textrm{degree} (q)=m$. Define \begin{displaymath} K(z) = z^{m-\ell} \frac{\widetilde{q}(z)/q(z)}{B(z)W(z)}. \end{displaymath} Then $|K|=1$ on \textbf{T} and $FK \ge 0,$ a.e. on \textbf{T}. Moreover, $K$ has poles of order 1 (or no pole) at each point $\beta_1...,\beta_n$ with respective residues $c_1,..,c_n$ and no other poles in $\Delta$ since any zeros of $W$ are 
also zeros of $\widetilde{q}$. Thus, \begin{displaymath} 1 = \int| F| d \sigma= \int FK d \sigma =\sum_{k=1}^n c_k F(\beta_k)\end{displaymath} whereas \begin{displaymath}  \textrm{Re   } [\sum_{k=1}^n c_k f(\beta_k)] =  \textrm{Re  } [\int fK d \sigma] \le 1\end{displaymath} for all $f \in H^1, ||f|| \le 1.$ Thus, $(F(\beta_1),...,F(\beta_n))$ is a boundary point of $\Lambda.$

\bigskip

The extreme and exposed points of $\Lambda$ are described in Theorems 2 and 4, respectively, below. We continue the notation established above.

\bigskip

\textbf{Theorem 2} A point $(F(\beta_1),...,F(\beta_n))$ in the boundary of $\Lambda$ is an extreme point if and only if $F$ is outer and so is of the form $F=q^2P$.

\medskip

\textbf{Proof.} Assume first that $F=q^2P$. Let $K$ be the kernel function for $F$. Since $F$ is outer it has no zeros in $\Delta$ and so $m$ must be the order of the zero of $K$ at the origin.  Thus, \begin{equation} K(z)B(z) =z^m \frac{\widetilde{q}(z)}{q(z)}. \end{equation}To show that $(F(\beta_1),...,F(\beta_n))$ is an extreme point of $\Lambda$ we assume that  \begin{equation} F(\beta_k)= \lambda F_1(\beta_k) + (1 -\lambda) F_2(\beta_k), k=1,...,n \end{equation} for some $\lambda, 0 < \lambda <1$ and some functions $F_1,F_2$ in $H^1$ of unit norm for which $(F_j(\beta_1),...,F_j(\beta_n)), j=1,2$ lies in the boundary of $\Lambda.$ We must show that\begin{displaymath}(F_1(\beta_1),...,F_1(\beta_n))=(F_2(\beta_1),...,F_2(\beta_n))\end{displaymath}  From (12) we have \begin{displaymath} F= \lambda F_1 + (1-\lambda)F_2 + Bh \end{displaymath} where $h \in H^{\infty}$. Multiply both sides by $K$. Note that since $K(0)=0$ the integral of $KBh$ over \textbf{T} with respect to Lebesgue measure is zero. This then gives \begin{eqnarray*} 1 & = &  \int |F| d \sigma =  \int FK  d\sigma \\ & = &  \lambda \int KF_1 d\sigma + (1-\lambda) \int K F_2 d\sigma \\  & \le & \lambda \int|F_1| d\sigma + (1 -\lambda)\int|F_2| d\sigma =1. \end{eqnarray*} Hence, \begin{equation} F_j K \ge 0 \textrm{  a.e. on \textbf{T}}, j=1,2. \end{equation} From (11) and (13) we have \begin{eqnarray*} 0 & \le & \frac{F_j KB}{B} =  \frac{\widetilde{q}(e^{it})}{q(e^{it})} \frac{e^{imt}F_j(e^{it})}{B(e^{it})}  \\ & = & \frac{e^{imt} F_j(e^{it})[\widetilde{q}(e^{it})/q(e^{it})] }{e^{int} P(e^{it})} [ \frac{e^{int} P(e^{it})}{B(e^{it})} ]\end{eqnarray*} However, the final term on the far right-hand side just above is positive on \textbf{T} by (6) and we deduce that \begin{displaymath} \frac{F_j(e^{it})[\widetilde{q}(e^{it})/q(e^{it})]}{P(e^{it}) e^{i(n-m)t}} \ge 0, \textrm{  on \textbf{T}},  j=1,2 .\end{displaymath}  Let $m_j$ be the order of the zero of $F_j$ at the origin and let $W_j$ denote the Blaschke factor of $F_j$ from (4). Thus there are polynomials $p_1,p_2$ of degrees $n-m-m_j$ with no zero in $\Delta$ such that \begin{equation} F_j(z)\frac{\widetilde{q}(z)}{q(z)}= p_j(z) \widetilde{p_j}(z)P(z)= \frac{\widetilde{p_j}(z)}{p_j(z)} p_j^2(z)P(z), j=1,2.\end{equation}   We may now equate the inner factors of the far sides of (14). This yields \begin{equation} z^{m_j}W_j(z) \frac{\widetilde{q}(z)}{q(z)}= \frac{\widetilde{p_j}(z)}{p_j(z)}. \end{equation}  (15) implies that $m_j=0$ and \begin{equation} \textrm{every zero of   }\widetilde{q} \textrm{  in  } \Delta  \textrm{  is a zero of the same or higher order of   } \widetilde{p_j}, j=1,2 \end{equation} Now substitute (14) into (12) and cancel $P(\beta_k)$ from both sides. It follows that the polynomial $W$ of degree $2(n-m)$ given by \begin{displaymath} W(z)= q(z) \widetilde{q}(z) - [\lambda p_1(z) \widetilde{p_1}(z) + (1-\lambda)p_2(z) \widetilde{p_2}(z)] \end{displaymath} vanishes at $\beta_k, k=1,..,n$. However, \begin{displaymath} z^{2(n-m)} \overline{W(1/\overline{z})} = W(z) \end{displaymath} so that $W$ also vanishes at the reflected points $1/\overline{\beta_k}, k=1,..,n$. Thus, $W$ has at least $2n$ zeros and so vanishes identically. Therefore,  \begin{displaymath} q(z) \widetilde{q}(z) = \lambda p_1(z) \widetilde{p_1}(z) + (1-\lambda)p_2(z) \widetilde{p_2}(z) \end{displaymath} Multiply both sides by $z^{m-n}$ and consider the result on \textbf{T}. We have \begin{equation}  |q(z)|^2  = \lambda |p_1(z)|^2 + (1-\lambda)|p_2(z)|^2, |z|=1. \end{equation} (17) implies that both $p_1$ and $p_2$ vanish at each zero of $q$ on \textbf{T} to the same or higher order. Thus, $p_1$ and $p_2$ have at least as many zeros in \textbf{T} as does $q$. When we combine this with (16) we find $q, p_1,p_2$ have the same zeros and so are the same polynomial up to a scale factor. This implies that $F_1=F_2$ and so $(F(\beta_1),...,F(\beta_n))$ is an extreme point of $\Lambda$.

\medskip

Conversely, suppose that  $(F(\beta_1),...,F(\beta_n))$ is an extreme point of $\Lambda$ and let $K$ be the kernel function for $F$. We must show that $F$ is outer. If it is not, then the deLeeuw-Rudin Theorem implies there are outer functions $F_1,F_2$ of $H^1$-norm 1 with \begin{equation} F = \frac{1}{2}F_1 + \frac{1}{2}F_2 .\end{equation}  We have \begin{eqnarray*} 2 & = & 2 \int |F| d \sigma = 2\int FK d\sigma \\ & = & \int F_1 K d \sigma + \int F_2 K d \sigma \\ & \le & \int |F_1| d \sigma = \int |F_2| d\sigma = 2. \end{eqnarray*} Hence, $ F_jK \ge 0, j=1,2$  a.e. on \textbf{T}. Now $(F(\beta_1),...,F(\beta_n))$ is an extreme point of $\Lambda$ and so (18) implies that \begin{equation} F(\beta_k) = F_1(\beta_k)=F_2(\beta_k), k=1,..,n. \end{equation} Therefore both $(F_1(\beta_1),...,F_1(\beta_n))$ and $(F_2(\beta_1),...,F_2(\beta_n))$ lie in the boundary of $\Lambda$. From Proposition 1 and the assumption that $F_1, F_2$ are outer, we learn that $F_j= q_j^2P, j=1,2$ where no zero of $q_1$ or $q_2$ lies in $\Delta.$. Suppose $q_j$ has $r_j \ge 0$ zeros on \textbf{T}; there is no loss of generality if we assume $r_1 \ge r_2$. Now if $\lambda_1,..,\lambda_r$ all lie in \textbf{T}, then  \begin{equation} \prod_1^r (e^{it} -\lambda_k)^2 = \mu e^{irt} \prod_1^r |e^{it}-\lambda_k|^2 \end{equation} where $\mu$ is a unimodular constant depending on $\lambda_1,..,\lambda_r.$
Therefore, the two inequalities $KF_j \ge 0, j=1,2$  imply that \begin{eqnarray*} 0  & \le & \frac{F_1K}{F_2K} = \frac{F_1}{F_2} = \frac{q_1^2}{q_2^2} \\ & = & \frac{V_1(e^{it})}{V_2(e^{it})} \frac{\prod_1^{r_1}(e^{it}- \lambda_k)^2}{\prod_1^{r_2}(e^{it} - \mu_k)^2}  \\ & = & C e^{i(r_1-r_2)t}\frac{V_1(e^{it})}{V_2(e^{it})}  \frac{\prod_1^{r_1}|e^{it}- \lambda_k|^2}{\prod_1^{r_2}|e^{it}- \mu_k|^2} \end{eqnarray*} where $V_1,V_2$ are polynomials all of whose zeros lie outside \textbf{T} and $C$ is a unimodular constant. Hence, \begin{displaymath} 0 \le C e^{i(r_1-r_2)t} \frac{V_1(e^{it})}{V_2(e^{it})}. \end{displaymath} However, \begin{displaymath}  Cz^{r_1-r_2} \frac{V_1(z)}{V_2(z)} \end{displaymath} is in $H^{\infty}$. Since it is non-negative on \textbf{T}, it must be constant. Hence, $r_1=r_2$ and $V_1=AV_2$ where $A$ is a constant.  That is, $q_1$ and $q_2$ have the same zeros outside \textbf{T}. Let $V$ be the polynomial with exactly these zeros so that $q_j=Vp_j, j=1,2$ where $p_1,p_2$ have all their zeros on \textbf{T}; let $r=r_1=r_2$ be the number of zeros of $q_1$ and $q_2$ on \textbf{T}. Refer now to (19). After cancelling the factor $PV^2$ that appears in both $q_1^2$ and $q_2^2$ we obtain \begin{equation} p_1^2(\beta_k)= p_2^2(\beta_k), k=1,..,n \end{equation} However,  all the zeros of $p_j, j=1,2$ lie on \textbf{T} and $0\le p_1^2/p_2^2$ on \textbf{T} so that
$\widetilde{p_j} = \mu p_j, j=1,2$ where $\mu$ is a unimodular constant. Thus \begin{displaymath}p_1^2(1/\overline{\beta_k})= p_2^2(1/\overline{\beta_k}), k=1,..,n \end{displaymath} which implies that $p_1^2=p_2^2$ at 2$n$ points and so $p_1^2 \equiv p_2^2$. This in turn implies that $F_1 = F_2$. This contradicts the assumption that $F$ is not outer. So $F$ must be outer and of the form $F=q^2P$ if $(F(\beta_1),...,F(\beta_n))$ is an extreme point of $\Lambda$. This completes the proof of Theorem 2.

 \bigskip
 
 The following shows that $\Lambda$ inherits a structure reminiscent of the unit ball of $H^1$; see [1; p. 125].

\bigskip

\textbf{Corollary 3} If $\vec{P}$ is in the boundary of $\Lambda$ and is not an extreme point, then there are extreme points $\vec{P_1}, \vec{P_2} \in \Lambda$ such that \begin{displaymath} \vec{P}= \frac{1}{2} \vec{P_1} + \frac{1}{2} \vec{P_2}.  \end{displaymath}

\medskip

\textbf{Proof} Let  $\vec{P}=(F(\beta_1),...,F(\beta_n))$ where $F \in H^1$ has $H^1$-norm 1. Since $\vec{P}$ is not an extreme point of $\Lambda$, we know from Proposition 1 and Theorem 2 that $F$ is not outer. The deLeeuw-Rudin Theorem implies there are outer functions $F_1,F_2$ of unit norm with \begin{displaymath} F= \frac{1}{2} F_1 + \frac{1}{2} F_2 \end{displaymath} Since $\vec{P}$ is in the boundary of $\Lambda$, it follows that $\vec{P_j}= (F_j(\beta_1),...,F_j(\beta_n)), j=1,2$ must also lie in the boundary of $\Lambda$.  Theorem 2  implies that $\vec{P_j}$ is an extreme point of $\Lambda$ since $F_j$ is outer.

\bigskip

\textbf{Theorem 4} An extreme point in the boundary of $\Lambda$ is an exposed point if and only if the polynomial $q$ from Theorem 2 has no zeros on \textbf{T}

\medskip

\textbf{Proof} Suppose $(F(\beta_1),...,F(\beta_n))$ is an exposed point of $\Lambda$. Then it is an extreme point so that $F=q^2P$ where $q$ has no zeros in $\Delta$.  Let $K$ be the kernel function for $F$ so that $FK \ge 0$ and $|K|=1$ on \textbf{T}.  Suppose that $q$ has a zero at $\lambda, |\lambda|=1$. Define \begin{displaymath} G(z) = A(-\lambda) \frac{z F(z)}{(z-\lambda)^2}= A(- \lambda)z (\frac{q(z)}{z-\lambda})^2P(z) .\end{displaymath} Then $G \in H^1$; we chose $A$ to be a positive scalar to make $||G||_1=1.$ Since \begin{displaymath} (-\lambda) \frac{e^{it}}{(e^{it}-\lambda)^2} \ge 0 \end{displaymath} we see that $KG \ge 0$ on \textbf{T}. Hence, \begin{displaymath} \sum_1^n c_k G(\beta_k) = \int KG d\sigma = \int |G| d \sigma =1. \end{displaymath} However, $(G(\beta_1,...,G(\beta_n))) \ne (F(\beta_1),...,F(\beta_n))$. This contradicts the fact that $(F(\beta_1),...,F(\beta_n))$ is an exposed point. Hence,  $q$ has no  zero on \textbf{T} 

\medskip

Conversely, suppose $F=q^2P$ where $q$ has no zero on \textbf{T}. Let $K$ be a kernel function for $F$ so that $FK \ge 0$ and $|K|=1$ on \textbf{T} and \begin{equation} 1 = \int KF d \sigma \ge \textrm{Re  } \int fK d \sigma \end{equation} for all $f \in H^1, ||f|| \le 1.$  Suppose $G \in H^1$ satisfies $||G||_1=1$ and  \begin{equation} \sum_1^n c_k G(\beta_k)= \int GK d \sigma =1. \end{equation}Then $(G(\beta_1),...,G(\beta_n))$ must lie in the boundary of $\Lambda$ and so Proposition 1 tells us that \begin{displaymath} G(z)= z^{\ell} W(z) p^2(z) P(z) \end{displaymath} where $p$ is a polynomial with no zeros in $\Delta.$ Moreover, (23) implies that  $GK \ge 0$ on \textbf{T}. Hence, $G/F \ge 0$ a.e. on \textbf{T}. Therefore,  for $|z|=1$, we have \begin{equation} 0 \le \frac{G(z)}{F(z)} = \frac{z^{\ell}W(z)p^2(z)P(z)}{q^2(z)P(z)} = z^{\ell}\frac{W(z)p^2(z)} {q^2(z)} \end{equation} But the last term on the right in (24) is in $H^{\infty}$ since $q$ has no zero on \textbf{T}. Thus the right-hand side of (24) is a bounded analytic function that is non-negative on \textbf{T} and so constant. This implies that $\ell =0$ and $W$ is a constant and, lastly, that $p^2 =A q^2$ for some constant $A$. Consequently $G=F$. Hence, $(F(\beta_1),...,F(\beta_n))$ is an exposed point of $\Lambda.$

\bigskip

\textbf{Corollary 5} A point $\vec{P}=(F(\beta_1),...,F(\beta_n))$ in the boundary of $\Lambda$ is an exposed point if and only if $F$ is the unique solution to the corresponding extremal problem: \begin{equation} \sup \{ \textrm{Re  } [\sum_1^n c_k f(\beta_k)]: f \in H^1, ||f||_1 \le 1 \}.  \end{equation}

\medskip

\textbf{Proof} Evidently if (25) has a unique solution $F$ in $H^1$, then $(F(\beta_1),...,F(\beta_n))$ is an exposed point of $\Lambda$. Conversely, suppose that $F_1$ and $F_2$ are two functions in the unit sphere of $H^1$ and $(F_1(\beta_1),...,F_2(\beta_n))=(F_2(\beta_1),...,F_2(\beta_n))$ is an exposed point of $\Lambda$. Theorem 4 implies that there are polynomials $q_1,q_2$ all of whose zeros lie outside \textbf{T} such that $F_j=q_j^2P, j=1,2$. Hence, \begin{equation} 0 \le \frac{F_1K}{F_2K} = \frac{F_1}{F_2} = \frac{q_1^2}{q_2^2}. \end{equation} However, the function on the far right-hand side of (26) is in $H^{\infty}$ since all the zeros of $q_2$ lie outside \textbf{T}. So it must be constant and therefore $F_1=F_2$.

\bigskip

\textbf{Further observation}  The assumption that the points $\beta_1,..,\beta_n$ are distinct is not really needed provided we follow the convention that if some $\beta_k$ is repeated, say, $r$ times, then the $r$ values $F(\beta_k)$ are understood to be \begin{displaymath} F(\beta_k), F'(\beta_k),...,F^{(r-1)}(\beta_k). \end{displaymath} With this understanding and much more complicated notation in the proofs, the conclusions of Theorems 2 and 3 hold.

\section*{References}

\begin{enumerate}

\item P. Duren, \textit{The Theory of $H^p$ Spaces}, Academic Press, New York, 1970

\item S. Ya. Havinson,  "Two papers on Extremal Problems in Complex Analysis," Amer. Math. Translation series 2, vol. 129, Providence, R.I., 1986

\item K. Hoffman, \textit{Banach Spaces of Analytic Functions}, Prentice-Hall, Englewood Cliffs, N.J., 1962

\end{enumerate}

\end{document}